\documentclass [12pt]{article}
\usepackage[T2A]{fontenc}
\usepackage[cp866]{inputenc}
\usepackage[english]{babel}
\begin{document}

\begin{center}
	{\bf I.V. Boykov, A.I. Boykova
		
		Analytical methods for solution of hypersingular and polyhypersingular integral equations
		}	
\end{center}
\begin{abstract}
We propose a method for transformating linear and nonlinear hypersingular integral equations into ordinary differential equations.  Linear and nonlinear polyhypersingular integral equations are transformed into partial differential equations. Well known that many types of differential equations can be solved in quadratures. So, we can receive analytical solutions for many types of linear and nonlinear hypersingular and	polyhypersingular integral equations.

\end{abstract}

MSC 45E05, 45G05

\section{Introduction}

Importance of solving hypersingular integral equations is justified by numerous applications and intense 
growth of the fields during the last century since Hilbert and Poincare created the theory of singular integral equations.
The theory is associated with numerous applications of singular and hypersingular integral equations, as well as with Riemann's boundary
value problem.  The Riemann's boundary value problem, singular, and hypersingular integral equations are broadly used as basic techniques of mathematical modeling in physics (quantum field theory \cite{Bogo}, theory of short and long-range interaction 
\cite{Br}, soliton theory \cite{Fad}), theory of elasticity and thermoelasticity \cite{Ia}, aerodynamics and electrodynamics 
\cite{Lif} and many other fields.

A closed-form solution of singular and hypersingular integral equations is only possible in
exceptional cases.  A comprehensive presentation and an extensive literature survey associated 
with all methods of solution of singular integral equations of the first and second kinds can be
found in \cite{Gakh, Mus,Ivan,Goh,Gol,Mich,Siz,Boy1,Boy2}.   
The methods of solution of hypersingular integral
equations are  less elaborated. In this paper we pay a special attention to solving linear
and nonlinear hypersingular integral equations of the first and second kinds as the equations of
this kind describe many problems in fluid dynamics and electrodynamics. 
For example, we refer to one-dimensional \cite{Cap} and multi-dimensional  Prandtl's equations for
a steady flow around a wing \cite{Lif,Osel}. In \cite{Gol1, Gol2,Lif1,Mand} Galerkin and collocation 
methods  for solving  hypersingular integral equations of the first and second kind are developed
and justified by imposing conditions on the kernels and the right-hand sides of the equations.

In \cite{Boy3} - \cite{Boy8} a collocation and spline-collocation  methods  for solving hypersingular and polyhypersingular
integral equations  with even singularity has been developed and justified, and
the nonlinear case is described in \cite{Boy6,Boy7,Boy8}.

A closed-form solution of  hypersingular integral equations is only possible in a few cases. In this direction we can remark the works \cite{Mart}, \cite{Mart1}, \cite{Mand}, in which the exact solution of the equation
 $$
 Hx = \frac{d}{dt}\left(\int\limits_{-1}^1\frac{x(\tau)}{\tau-t}d\tau\right) = f(t), \, -1 < t < 1,
 \eqno (1.1)
 $$
 under conditions $x(\pm 1)=0$, is given. 
 
 In the book \cite{Mand} the exact solution of the equation 
 $$
 x(t) - \frac{\alpha}{\pi}(1-t)^{1/2} \int\limits^1_{-1}
 \frac{x(\tau)}{(\tau-t)^2}d\tau = f(t), \, -1 < t < 1, \  \alpha=const,
 \eqno (1.2)
 $$
 under conditions $x(\pm 1)=0$, is given. 
 
{\it Remark.} Unlike hypersingular integral equations, in the theory of singular integral equations, wide classes of equations are distinguished for which solutions are constructed in closed form. In addition to the characteristic singular integral equations, solutions in closed form are constructed for a series of complete singular integral equations \cite{Sam}, \cite{Gakh}.

In this paper, we propose an approach to solving linear and nonlinear  hypersingular integral
equations and polyhypersingular integral equations. We investigate conditions under which singular and hypersingular integral equations can be transformed into differential equations. It  allows us to use the big arsenal of methods for solution of differential equations to solving singular and hypersingular integral equations.

The paper is constructed in following way. In item two some definitions, used in the paper, are given. In item three a method for transformation hypersingular integral equations into ordinary differential equations is given. In item four considered  method for transformation nonlinear hypersingular integral equations into nonlinear ordinary differential equations.
In item five considered  method for transformation  hypersingular integro-differential   equations into  ordinary differential equations.
In item six investigated method for transformation
 polyhypersingular integral equations into partial differential equations.

\section{ Definitions}

Let us recall definitions of hypersingular integrals.
Hadamard \cite{Had} introduced a new type of  integrals:

{\bf Definition 2.1} The integral of the type  
$$
	\int\limits_a^b \frac{A(x) \, dx}{(b-x)^{p+\alpha}}
	\eqno (2.1)
$$
	for an integer $p$ and $0<\alpha<1$ one defines a value 
	of the above integral (''finite part'') as the limit of the sum
	$$ \int\limits_a^x \frac{A(t) \, dt}{(b-t)^{p+\alpha}} +
	\frac{B(x)}{(b-x)^{p+\alpha-1}}, $$
	as $x\to b$ assuming that $A(x)$ has $p$ derivatives in the neighborhood 
	of point  $b$. Here $B(x)$ is any function that satisfies the following two conditions:

		i) The above limit exists;
		
	ii) $B(x)$ has at least $p$ derivatives in the neighborhood of a point $x=b$.

An arbitrary choice of $B(x)$ is unaffected by the value of the limit 
in the condition (i).
The condition (ii) defines values
of the $(p-1)$ first derivatives of $B(x)$ at a point $b$,
so that an arbitrary additional term in a numerator is an infinitely small
quantity, at least of order $(b - x)^p$.

Chikin \cite{Chi} introduced the definition of the Cauchy - Hadamard type 
integral that generalizes the notion of the singular
integral in the Cauchy principal sense and in the Hadamard sense.

{\bf Definition 2.2.}
	The Cauchy - Hadamard principal sense of the integral
	$$ \int\limits_a^b \frac{\varphi(\tau) \, d\tau}{(\tau-c)^p}, \quad
	a<c<b, $$
	is defined as the limit of the following expression
	$$ \int\limits_a^b \frac{\varphi(\tau)\, d\tau}{ (\tau-c)^p} = \lim_{v\to
		0} \left[ \int\limits_a^{c-v} \frac{\varphi(\tau) \, d\tau}{(\tau-c)^p}
	+ \int\limits_{c+v}^b \frac{\varphi(\tau)\, d\tau}{(\tau-c)^p}
	+ \frac{\xi(v)}{v^{p-1}}\right], $$
	where $\xi(v)$ is a function chosen in such a way as to provide an existence of the above limit.

{\bf Definition 2.3.} The integral of the type  
	$ \int\limits_a^b \frac{\varphi(\tau)\, d\tau}{ (b-\tau)^p}$ 
	is defined by
	$$
	\int\limits_a^b \frac{\varphi(\tau)\, d\tau}{ (b-\tau)^p}=\lim_{v\to 0}
	\left[ \int\limits_a^{b-v} \frac{\varphi(\tau) \, d\tau}{(\tau-t)^p}
	+\frac{\xi(v)}{v^{p-1}}+\xi_1(v) \ln v \right],
	$$
	where $\xi(v)$ has $p-1$ derivatives, and  the $(p-1)$th derivative is continuous 
	in the neighborhood of zero, and $\xi_1(v)$ satisfies the  Dini-Lipschitz condition.
	One can choose functions $\xi(v)$ and $\xi_1(v)$ in such a way as to provide an existence of the above limit.

In this section we will list several classes of functions, which will be used later.

Let $\gamma$ be the unit circle: $\gamma = \{z:|z|=1\}.$

To measure the continuity of a function $f \in C[a,b]$ we proceed as follows \cite{Lor}. We consider the first difference with step $h$ 
$$
\Delta_h(f(x)) = f(x+h) - f(x)
$$
of the function $f$ and put
$$
\omega(f,\delta) = \omega(\delta) = \max\limits_{x,h  (|h| \leq \delta)}|f(x+h) - f(x)|.
$$
The function $\omega(\delta),$ called the modulus of continuity of $f,$ is defined for $0 \leq \delta \leq b-a.$

{\bf Definition 2.4.} A function $f,$ defined on $\Delta = [a,b]$ or $\Delta = \gamma,$ satisfies a Lipschits condition with constant $M$ 
and exponent $\alpha,$ or belongs to the class $H_\alpha(M),$ $M \geq 0,$ $0 < \alpha \leq 1,$ if 
$$
|f(x') - f(x'')| \leq M|x'-x''|^\alpha, \, x', x'' \in \Delta.
$$

{\bf Definition 2.5.} The class $W^r(M,\Delta),$ $r=1,2,\ldots,$ $\Delta = [a,b]$ or $\Delta = \gamma,$ consists of all functions $f \in C(\Delta),$ 
which have an absolutely continuous derivative $f^{(r-1)}(x)$ and piecewise derivative $f^{(r)}(x)$ with $|f^{(r)}(x)| \le M.$

{\bf Definition 2.6. } \cite{Lor}. \ Let  $r=0,1,\ldots,$  $M_i \geq 0,$  $i=0,1,\ldots,\\ r+1,$ let $\omega$ be a modulus of continuity and let $\Delta=[a,b]$ 
or $\Delta= \gamma.$ Then $W^r_\omega = W^r_\omega(M_0,\ldots, M_{r+1}; \Delta)$ is the set of all functions $f \in C(\Delta),$ which have continuous
derivatives $f, f',\ldots, f^{(r)}$ on $\Delta,$ satisfying
$$
|f^{(i)}(x)| \leq M_i, \, x \in \Delta, \, i=0,1,\ldots,r, \, \omega(f^{(r)}, \delta) \leq M_{r+1}\omega(\delta).
$$
They write $W^r H_\alpha$ if $\omega(\delta) = \delta^\alpha,$ $0 < \alpha \leq 1.$

Let us consider functions $f(x_1,\ldots,x_l)$ of $l$ variables on $\Delta,$ where $\Delta$ is either an $l$-dimensional parallelepiped (that is, the product of $l$
intervals $a_k \leq x_k \leq b_k,$ $k=1,2,\ldots,l)$ or an $l$-dimensional torus $\gamma^{(l)},$ the product of $l$ circles $\gamma.$ The modulus of continuity 
$\omega(f, \delta)$ are defined as the maximum of 
$$|f(y_1,\ldots,y_l) - f(x_1,\ldots,x_l)| \ {\rm for \ }
|y_k-x_k| \leq \delta, k=1,2,\ldots,l.
\eqno (2.1)
$$

In \cite{Lor} was defined classes $W^r_{l \omega} = W^r_{l \omega}(M_0,\ldots, M_{r+1}; \Delta)$ of functions $f(x_1,\ldots,x_l)$ of $l$ variables
on $l$-dimensional set $\Delta,$ which is either a parallelepiped or a torus. 

{\bf Definition 2.7.} \cite{Lor}. A function $f$ belongs to $W^r_{l \omega}(M_0,\ldots, M_{r+1}; \Delta)$ if and only if all its partial derivatives $D^jf$ of 
order $j=0,1,\ldots,r$ exist and continuous and satisfy the following conditions:
For each partial derivative $D^j f$ of order $j,$ $\|D^j f\| \leq M_j,$ $j=0,1,\ldots,r,$ and in addition for each derivative of order 
$r,$ $\omega(D^r f, \delta) \leq M_{r+1} \omega(\delta).$

They write $W^r_l H_\alpha$ if $\omega(\delta) = \delta^\alpha,$ $0 < \alpha \leq 1.$

If coefficients $A$ and $M$ are not essential we use designations
$H_{\alpha},$ \ \  $W^r H_{\alpha},$ \ \ $W^{r}_l H_{\alpha}$ \ \ instead of \
$H_{\alpha}(M,\Delta),$\ \ $W^r H_{\alpha}(M,\Delta),$ \\ $W^{r}_l H_{\alpha }(M,\Delta)$
respectively.

Similarly, definitions are introduced in the  case when an arbitrary smooth bounded contour is taken as the contour of $ \gamma $.

\section{ Hypersingular integral equations }

Let us consider 
hypersingular integral equation
$$
a(t)x(t) + \frac{1}{\pi i}\int_L \frac{h(t,\tau)x(\tau)d\tau}{(\tau-t)^p} +
\int_L k(t,\tau) x(\tau)d\tau = f(t), \, p=2,3,\ldots
\eqno (3.1)
$$
where $L$ is a smooth closed contour in the plane of a complex variable.

We denote by $ D^+ (D^-) $ the inner (outer) domain with respect to the contour $ L. $
We denote by $ \bar D^+ $ the closure of the domain $ D^+. $ Denote by $ G $ an open domain such that $ \bar D^+ \subset G. $

At first we  consider the singular integral equation
$$
a(t)x(t) + \frac{b(t)}{\pi i}\int_L \frac{x(\tau)d\tau}{(\tau-t)^p} = f(t).
\eqno (3.2)
$$

We will seek a solution of the equation (3.2) in the class of functions $x(z), z \in \bar D^+ $, which are analytical in $D^+$ and satisfied the condition:
 $ x(t) \in W^{p-1}H_\alpha (M,L).$
 
 Let the following conditions be satisfied for functions $a(t), b(t), f(t)$:\\
1) functions $a(z),$ $b(z),$ $f(z)$ are analytical in $G;$\\
2) $a(z)\ne 0,$   $z \in D^+.$

We will prove that, under these conditions, the equation (3.2) has at least 
 $p-1$ linearly independent solutions.

Assume that the equation  (3.2) has a solution $  x^*(t) \in W^{p-1}H_\alpha(M,L).$

As  $a(z) \ne 0,$
 $z \in G,$ that 
$$
x^*(t) = \frac{f(t)}{a(t)} - \frac{b(t)}{a(t)} \frac{1}{\pi i}\int_L \frac{x^*(\tau)d\tau}{(\tau-t)^p}d\tau.
$$

From definition of hypersingular integral
 \cite{Chi} , we have 
$$
x^*(t) = \frac{f(t)}{a(t)} - \frac{b(t)}{a(t)} \frac{1}{\pi i(p-1)!}\int_L \frac{x^{*(p-1)}(\tau)d\tau}{\tau-t}d\tau.
$$

Function
$$
\Psi(z) = \frac{1}{\pi i}\int_L \frac{x^*(\tau)d\tau}{(\tau-z)^p}
=\frac{1}{\pi i(p-1)!}\int_L \frac{x^{*(p-1)}(\tau)d\tau}{\tau-z}
$$
is analytical for $z \not \in L$.

So, the function
$$
x^*(z) = \frac{f(z)}{a(z)} -\frac{1}{(p-1)!} \frac{b(z)}{a(z)} \frac{1}{\pi i} 
\int\limits_{L} \frac{x^{*(p-1)}(\tau)d\tau}{\tau-z},
$$
is analytical in $D^+.$

Consequently, the function $ x^* (t) $ is the boundary value of the analytic function
$$
\frac{f(z)}{a(z)} - \frac{b(z)}{a(z)} \frac{1}{\pi i}\int_L \frac{x^*(\tau)d\tau}{(\tau-z)^p}d\tau=
$$
$$
= \frac{f(z)}{a(z)} - \frac{b(z)}{a(z)} \frac{1}{\pi i(p-1)!}
\int_L \frac{x^{*(p-1)}(\tau)d\tau}{\tau-z}d\tau.
$$

It can be extended to the domain $ D^+ $ by the function
$$
x^*(z) = \frac{f(z)}{a(z)} - \frac{b(z)}{a(z)} \frac{1}{\pi i}\int_L \frac{x^*(\tau)d\tau}{(\tau-z)^p}d\tau.
$$

The function  $x^*(z)$ is analytical in the domain  $D^+$. In addition,  $x^*(t) \in W^{p-1}H_\alpha(M,L)$ by assumption.

Then \cite{Chi}
$$
\frac{1}{\pi i}\int_L \frac{x^*(\tau)d\tau}{(\tau-t)^p}d\tau
= \frac{1}{\pi i(p-1)!} \frac{d^{p-1}}{d t^{p-1}}\int_L \frac{x^*(\tau)d\tau}{\tau-t}d\tau. 
$$

Well  known \cite{Gakh},  that for functions, analytical in $D^+$ and satisfying the condition
 $x^*(t) \in H_\alpha$,  the equality 
$$
\frac{1}{\pi i} \int_L \frac{x^*(\tau)}{\tau-t}d\tau = x^*(t)
$$
is valid.

So, we have
$$
\frac{1}{\pi i} \int_L \frac{x^*(\tau)}{(\tau-t)^p}d\tau =\frac{1}{(p-1)!}
\frac{d^{p-1}}{d t^{p-1}}x^*(t).
$$

Thus, the solution $x^*(t)$ of the equation (3.2) satisfies the differential equation
$$
\frac{b(t)}{(p-1)!} \frac{d^{p-1}}{d t^{p-1}} x^*(t) + a(t) x^*(t) = f(t).
\eqno (3.3)
$$

Thus, if the equation (3.2) has a solution $x(t) \in W^{p-1}H_\alpha(M,L)$ and conditions 1) - 2) are satisfied, then this solution satisfies the differential equation
$$
\frac{d^{p-1} x(t)}{d t^{p-1}} + (p-1)! \frac{a(t)}{b(t)}x(t) =
\frac{(p-1)! f(t)}{b(t)}.
\eqno (3.4)
$$

This is ordinary differential equation, which has 
 $p-1$ linearly independent solutions.
 
Now we will find solutions of the equation (3.4) for different  values of   $p.$

Let $p=1.$ Then equation (3.4) is transformed to the following functional equation
$$
x(t) \left(1+\frac{a(t)}{b(t)}\right) = \frac{f(t)}{b(t)}.
$$

This equation has the solution 
$$
x(t) = \frac{f(t)}{a(t) + b(t)}.
$$

If $a(z) + b(z) \ne 0,$ $z \in G,$  the function $x(t)$ is analytical in the domain $G.$ This function  is a   solution of the singular equation  (3.2). In that  we can directly verify by substituting $ x (t) $ in (3.2). 

Let $p=2.$ Hypersingular integral equations with the second-order singularity find wide applications in physics and engineering.

Hypersingular integral equation
$$
a(t)x(t) + \frac{b(t)}{\pi i} \int_L \frac{x(\tau)}{(\tau-t)^2}d\tau = f(t)
\eqno (3.5)
$$
is transformed to  ordinary differential equation
$$
\frac{dx(t)}{dt} + \frac{a(t)}{b(t)}x(t) = \frac{f(t)}{b(t)}.
$$

The last equation has the solution
$$
x(t) = e^{-\int\limits^t_{t_0} p(\tau) d\tau}
[x(t_0)+ \int\limits^t_{t_0}q(\tau) e^{\int\limits^\tau_{t_0} p(\tau) d\tau}d\tau],
\eqno (3.6)
$$
where $p(t) = \frac{a(t)}{b(t)},$ $q(t) = \frac{f(t)}{b(t)};$ integrals are calculated on closed smooth contour  $L,$ $t_0 \in L.$

It follows from (3.6) that in order to obtain a single-valued solution of (3.5) it is required to give an initial condition at some point $ t_0 \in L. $ 

{\bf Example.}

Let us consider the equation
$$
x(t) + \frac{2+t^2}{\pi i} \int\limits_L \frac{x(\tau)}{(\tau-t)^2}d\tau = 
(5+2t^2) e^{2t},
\eqno (3.7)
$$
where $L$ is a closed smooth contour, in the exterior of which there are points $\pm i, \pm \sqrt{2} i.$

The exact solution of this equation is
$$
x(t) = e^{2t}.
$$

The equation 
 (3.7) is transformed to differential equation
$$
\frac{dx(t)}{dt} + \frac{1}{2+t^2}x(t) = \frac{5+2t^2}{2+t^2}e^{2t}.
\eqno (3.8)
$$

A particular solution of this equation is the function $x_{pr} = e^{2t}.$

The general solution of the homogeneous equation
$$
\frac{dx(t)}{dt} + \frac{1}{2+t^2}x(t) = 0
$$
is the function $x_c(t) = A \rm{exp} \{-\frac{1}{\sqrt{2}} arctg \frac{t}{\sqrt{2}}\},$ where $A=\rm{const.}$ 

The function $arctg z$ is expressed in terms of the natural logarithm by the formula 
$$
arctg z = -\frac{i}{2} \rm{Ln}\left(\frac{1+iz}{1-iz}\right).
$$

So, the function  $arctg z$ is analytic in the plane of a complex variable everywhere, with the exception of points $\pm i.$

We denote by $ G $ the bounded closed domain in whose exterior the points $\pm i, \pm \sqrt{2} i$ lie.

In this domain the functions
$x^*(t) = e^{2t} + A e^{-\frac{1}{\sqrt{2}} arctg\frac{t}{\sqrt{2}}}$ and $\frac{1}{1+t^2}$ are analytic.

So, the function $x^*(t)$ is a solution of hypersingular integral equation  (3.7).

A direct check confirms this assertion.

The proposed method is extended to compound hypersingular integral equations and nonlinear hypersingular integral equations.

As an example, consider the composite hypersingular integral equation
$$
\frac{a_p(t)}{\pi i} \int\limits_L \frac{x(\tau)d\tau}{(\tau-t)^p} + 
\frac{a_{p-1}(t)}{\pi i}\int\limits_L \frac{x(\tau)d\tau}{(\tau-t)^{p-1}} + \cdots +
$$
$$
+ \frac{a_1(t)}{\pi i} \int\limits_L \frac{x(\tau)d\tau}{\tau-t} + a_0(t) x(t) = f(t),
\eqno (3.9)
$$
where $L$ is a smooth closed contour.

Let the contour $L$  be located in a bounded closed domain $ G. $
Suppose that the equation (3.9) has a solution $ x^* (t) $ that is analytic in the domain $ G. $

Then, as follows from the definition of  hypersingular integrals, the solution $ x^* (t) $ of equation (3.9) is also the solution of the differential equation
$$
a_p(t)\frac{d^{p-1}x(t)}{dt^{p-1}} + a_{p-1}(t)\frac{d^{p-2}x(t)}{dt^{p-2}} + \cdots +
$$
$$
+ a_1(t)\frac{dx(t)}{dt} + a_0(t) x(t) = f(t).
\eqno (3.10)
$$

If the equation (3.10) has a solution that is representable in analytic form and is analytic in the domain $ G $, then this solution is also a solution of the hypersingular integral equation (3.9).

In particular, equations with constant coefficients and with an analytic right-hand side possess this property.

Let us return to the equation (3.1).

We will consider the equation 
$$
a(t)x(t) + \frac{1}{\pi i}\int_L \frac{h(t,\tau)x(\tau)d\tau}{(\tau-t)^p} = f(t).
\eqno (3.11)
$$

The equation (3.1) is investigated with similar way.

We will show that equation (3.11), under the conditions listed above, has at least  $ p-1 $ linearly independent solutions.

Using the definition of hypersingular integrals \cite{Chi}, equation (3.11) can be represented as

$$
a(t)x(t) +
$$
$$
+ \frac{1}{(p-1)!\pi}\int\limits_L a \frac{h_{0,p-1}(t,\tau)x(\tau)+h_{0,p-2}(t,\tau)x'(\tau)+\ldots+ h_{0,0}(t,\tau)x^{(p-1)}(\tau)d\tau}{\tau-t} =
$$
$$
= f(t),
\eqno (3.12)
$$ 
where $h_{0,i}(t,\tau)=\frac{\delta^{i}h(t,\tau)}{\delta\tau^i}, i=0,1,\ldots,p-1.$

From the theorem on limit values of singular integrals \cite {Gakh},  follows the equality
$$
\frac{1}{\pi i}\int_L \frac{h_{0,i}(t,\tau)x^{*(p-1-i)}(\tau)d\tau}{\tau-t} =h_{0,i}(t,t) x^{*(p-1-i)}(t).
\eqno (3.13)
$$

In proving equality (3.13), it suffices to restrict ourselves to considering the integral

$$
\frac{1}{\pi i} \int\limits_L \frac{h(t,\tau)x(\tau)}{\tau-t}d\tau,
$$
where the function $ h (t, \tau) $ (in both variables) is analytic in the open domain $ G, $ in which the contour of $ L$ is located. 

Obviously,
$$
\frac{1}{\pi i} \int\limits_L \frac{h(t,\tau)x(\tau)}{\tau-t}d\tau = 
$$
$$
= h(t,t) x(t) + \frac{1}{\pi i} \int\limits_L 
\frac{h(t,\tau) x(\tau) - h(t,t) x(t)}{\tau-t}d\tau.
$$

The function $ \psi (t, z) = \frac {h (t, z) x (z) - h (t, t) x (t)} {z-t} $
analytical inside the contour  $ L. $

Let   $ R (t, \rho) $ be the  circle with  sufficiently small radius $ \rho $ and with the center in the point $ t $.
Denote by $ L_ {\rho} $ the contour consisting of the part of the circle
$ R (t, \rho), $ located inside of the contour $ L $, and part of the contour  $ L, $ remaining after removing the part of $ L$ that  fell into the circle
$ B (t, \rho). $ We denote the latter part as $ L^*_\rho. $

Easy to see that
$$
\frac{1}{\pi i} \int\limits_{L_\rho} \psi(t,\tau)d\tau = 0
$$
and
$$
\frac{1}{\pi i} \int\limits_{L} \frac{h(t,\tau) x(\tau) - h(t,t) x(t)}{\tau-t}d\tau =
\lim\limits_{\rho \to 0} \left[\frac{1}{\pi i} \int\limits_{L^*_\rho} \frac{h(t,\tau) x(\tau) - h(t,t) x(t)}{\tau-t}d\tau \right] =
$$
$$
= \lim\limits_{\rho \to 0} \left[\frac{1}{\pi i} \int\limits_{L_\rho} \frac{h(t,\tau) x(\tau) - h(t,t) x(t)}{\tau-t}d\tau  - \frac{1}{\pi i} \int\limits_{R^+(t, \rho)} \frac{h(t,\tau) x(\tau) - h(t,t) x(t)}{\tau-t}d\tau \right] =
$$
$$
= \lim\limits_{\rho \to 0} \left[-\frac{1}{\pi i} \int\limits_{R^+(t, \rho)} \frac{h(t,\tau) x(\tau) - h(t,t) x(t)}{\tau-t}d\tau \right]= 0.
$$

Here $ R^+ (t, \rho)$ is the  part of the circle $ R (t, \rho), $ having a non-empty intersection with the domain $ \bar D^+. $ 

So,
$$
\frac{1}{\pi i} \int\limits_L \frac{h(t,\tau)x(\tau)}{\tau-t}d\tau=h(t,\tau)x(\tau).
$$

Using the  equality (3.13), we transform the hypersingular integral equation (3.11) into an ordinary differential equation
$$
\frac{b(t)}{(p-1)!}\left[\sum\limits_{i=0}^{p-1}h_{0,i}(t,t)x^{(p-1-i)}(t)\right]  + a(t) x(t) = f(t).
\eqno (3.14)
$$

Thus, the solutions of equation (3.11) belonging on the contour $L$ to the class of functions $ W^{p-1} H_\alpha, $ are contained among the solutions of the differential equation (3.14).

The reverse is also true. Let the differential equation (3.14) has in the domain $ G $, in which the contour $ L$ is located,  an analytical solution $ x^* (t). $ Then the function $ x^{*(p-1)} (t) $ is analytic in $ D^+ $ and continuously extends to the contour $ L. $ Under these conditions, equality (3.13) holds. Substituting this equality into equation (3.14), we obtain the  equation (3.12), which is equivalent to the equation (3.11).

\section{Nonlinear  hypersingular integral equations }

Let us consider  non-linear hypersingular integral equations as
$$
\pounds\left(\frac{1}{\pi i}\int\limits_L
\frac{x(\tau)d\tau}{(\tau-t)^p}, \frac{1}{\pi i}\int\limits_L \frac{x(\tau)d\tau}{(\tau-t)^{p-1}}, \cdots,\frac{1}{\pi i}\int\limits_L \frac{x(\tau)d\tau}{\tau-t}, x(\tau)\right) = f(t),
\eqno (4.1)
$$
where $\pounds $ is non-linear operator, $L$ is a smooth closed contour.

Let $ G  $ be a bounded closed domain in which the contour $ L $ is located.

If the equation (4.1) has an analytic solution in $ G $, then this solution satisfies the differential equation
$$
\pounds_1(x^{(p-1)}(t), x^{(p-2)}(t), \cdots, x'(t), x(t)) = f(t),
\eqno (4.2)
$$
where $\pounds_1$ is the operator obtained from the operator $ \pounds$ by replacing the hypersingular operators and the singular operator into the corresponding differential operators.

If the solution $x^*(t)$ of the equation (4.2) is analytic in $ G $, then this solution  is also the solution of equation (4.1).

{\bf Example.}

Let us consider the equation
$$
\frac{1}{\pi i}\int\limits_\gamma
\frac{x(\tau)d\tau}{(\tau-t)^2} -\left( \frac{1}{\pi i}\int\limits_\gamma \frac{x(\tau)d\tau}{\tau-t}\right)^2 - 3\frac{1}{\pi i}\int\limits_\gamma \frac{x(\tau)d\tau}{\tau-t} = -4, \, 
t \in \gamma,
\eqno (4.3)
$$
where $\gamma$ is a bounded closed smooth curve.
 
The equation (4.3) has the solution
$$
x(t) = \frac{c_1 - 4c_2 e^{5t}}{c_1+c_2 e^{5t}}, \, t \in \gamma.
\eqno (4.4)
$$

Using the method proposed above, we reduce the equation (4.3) to the differential equation
$$
\frac{dx(t)}{dt} - x^2(t) - 3x(t) = -4, \, t \in \gamma.
\eqno (4.5)
$$

It is known \cite{Kam},
that the solution of the differential equation
$$
y'(t) - y^2(t) - 3y(t) + 4 =0, \, t \in R
$$
is the function
$$
y=\frac{c_1-4c_2 e^{5t}}{c_1+c_2 e^{5t}}, \, t \in R.
$$

It can be established by direct verification that the function
$$
x(z)=\frac{c_1-4c_2 e^{5z}}{c_1+c_2 e^{5z}}
$$
is the solution of the equation 
$$
x'(z) - (x(z))^2 - 3x(z) + 4 = 0
$$
at any point $ z $ of the plane of the complex variable $ z. $

The function
$$
x(z)=\frac{c_1-4c_2 e^{5z}}{c_1+c_2 e^{5z}} = \frac{c-4 e^{5z}}{c+e^{5z}}
$$
is analytic in the plane of the complex variable everywhere except for the point $ z = \frac{1}{5} \rm {Ln} (-c), c=\rm {const}. $
 
We denote by  $\rm{Im_G}\varphi (z)$  the set of values of the function $\varphi(z)=-e^{5z}$ for $ z \in G. $

Thus, the function $x(t) = \frac{c-4e^{5t}}{c+e^{5t}},$
$t \in \gamma$  is a solution of equation (4.3) for  $c  \notin \rm{Im_G}\varphi (z),$.

\section{ Hypersingular integro-differential equations }

Let us consider  hypersingular integro-differential equation
$$
\sum\limits^n_{k=0} a_k(t) x^{(k)}(t) + \sum\limits^m_{l=0} \frac{a_l(t)}{\pi i}
\int\limits_L \frac{x^{(l)}(\tau)}{(\tau-t)^p}d\tau = f(t)
\eqno (5.1)
$$
with initial conditions
$$
x^{(v)}(c) = x_v, \, v=0,1,\ldots,r-1, \, r=\rm{max}(n,m),
\eqno (5.2)
$$
where  $c $ is a point on a smooth closed bounded contour $L.$

We assume that the functions $a_k(t),$ $k=0,1,\ldots,n,$ $b_k(t),$ $k=0,1,\ldots,m,$ $f(t)$ are analytic in the open domain $ G $ containing the contour $L.$

The solution of the Cauchy problem (5.1), (5.2) we will be sought in the class of functions analytic in the domain $ G. $

Using the definition of the hypersingular integral, we reduce equation (5.1) to the form
$$
\sum\limits^n_{k=0} a_k(t) x^{(k)}(t) + \sum\limits^m_{l=0} \frac{b_l(t)}{\pi i}
\frac{1}{(p-1)!}\int\limits_L \frac{x^{(l+p-1)}(\tau)}{\tau-t}d\tau = f(t).
\eqno (5.3)
$$

We assume that functions
 $x^{(k)}(t),$ $k=0,1,\ldots,r+p-2,$ 
are analytic in the domain  $G.$ Under these conditions the hypersingular integro-differential equation (5.3) is reduced to differential  equation 
$$
\sum\limits^n_{k=0} a_k(t) x^{(k)}(t) + \sum\limits^m_{l=0} b_l(t)
x^{(l+p-1)}(t) = f(t).
\eqno (5.4)
$$

Order of differential equation  (5.4) is equal to $s=\rm{max}(n, m+p-1).$

To solve it, a $r$ initial conditions are given.

Thus, the Cauchy problem (5.4), (5.2) is solved under $ s-r $ free parameters.

\section{ Polyhypersingular integral equations }

In this section classes of polyhypersingular integral equations, for which it is possible to obtain a solution in the analytical form, are distinguished. 

Below, for simplicity of notation, we restrict ourselves to  consideration of  bihypersingular integral equations.

We will consider the bihypersingular integral equation
$$
a(t_1,t_2) x(t_1,t_2) + \frac{b(t_1,t_2)}{\pi i} 
\int\limits_{\gamma_1} \frac{x(\tau_1, t_2)}{(\tau_1-t_1)^p}d\tau_1 +
$$
$$
+ \frac{c(t_1,t_2)}{\pi i} 
\int\limits_{\gamma_2} \frac{x(t_1, \tau_2)}{(\tau_2-t_2)^p}d\tau_2 -
\frac{d(t_1,t_2)}{\pi^2} 
\int\limits_{\gamma_1} \int\limits_{\gamma_2} \frac{x(\tau_1, \tau_2)d\tau_1 d\tau_2}{(\tau_1-t_1)^p (\tau_2-t_2)^p} = 
$$
$$
= f(t_1, t_2), \, (t_1, t_2) \in (\gamma_1 \times \gamma_2),
\eqno (6.1)
$$
where $\gamma_i $ is a smooth closed contour in the plane of a complex variable $z_i,$ $i=1,2.$ Let $\gamma = \gamma_1 \times \gamma_2.$

We denote by $ G_i $ a closed bounded domain in the plane $ z_i $ such that the contour $ \gamma_i $ lies inside $ G_i, $ $ i = 1,2. $

Let $G=G_1 \times G_2.$ Let functions $a(t_1,t_2), b(t_1,t_2), c(t_1,t_2), d(t_1,t_2), f(t_1, t_2)$ are analytic in the domain $G.$

Suppose that equation (6.1) has a solution $x^*(t_1, t_2)$ that is analytic in the domain $ G $. Then the function $x^*(t_1, t_2)$ is a solution of the differential equation
$$
a(t_1,t_2) x(t_1,t_2) + \frac{b(t_1,t_2)}{(p-1)!} 
 \frac{\partial^{p-1}x(t_1, t_2)}{\partial t_1^{p-1}}
+ \frac{c(t_1,t_2)}{(p-1!)} 
 \frac{\partial^{p-1}x(t_1, t_2)}{\partial t_2^{p-1}} -
$$
$$
-
\frac{d(t_1,t_2)}{((p-1)!)^2} 
 \frac{\partial^{2p-2} x(t_1, t_2)}{\partial t_1^{p-1} \partial t_2^{p-1}}
= f(t_1, t_2), \, (t_1, t_2) \in (\gamma_1 \times \gamma_2).
\eqno (6.2)
$$

Now we will use analogs of the Sokhotsky - Plemel formulas for multiple integrals of Cauchy type.

Let the contour $ \gamma_i $ divide the plane of the complex variable $ z_i, i = 1,2, $ into two parts: the inner $D^+_i$ and the outer $D^-_i.$ Then $\gamma = \gamma_1 \times \gamma_2$  is the boundary of regular bicylindrical domains $D^{\pm \pm} = D_1^{\pm} \times D_2^{\pm}.$

Consider the double Cauchy-type integral
$$
\Phi(z_1, z_2) = \frac{1}{(2\pi i)^2}
\int\limits_\gamma \frac{\varphi(\tau_1,\tau_2)d\tau_1 d\tau_2}{(\tau_1 - z_1)(\tau_2 - z_2)}.
\eqno (6.3)
$$

We will use  designations
$$
S_1 \varphi = \frac{1}{\pi i} \int\limits_{\gamma_1}
\frac{\varphi(\tau_1,t_2)d\tau_1}{\tau_1 - t_1},
$$
$$
S_2 \varphi = \frac{1}{\pi i} \int\limits_{\gamma_2} 
\frac{\varphi(t_1,\tau_2)d\tau_2}{\tau_2 - t_2},
$$
$$
S_{12}\varphi = -\frac{1}{\pi^2} \int\limits_\gamma
\frac{\varphi(\tau_1, \tau_2)d\tau_1 d\tau_2}{(\tau_1-t_1)(\tau_2 - t_2)}.
$$

We denote by $\Phi^{\pm \pm}(t_1, t_2)$ the limit values of the integral (6.3) when the point $(z_1, z_2) \in D^{\pm \pm}$ tends to the point $(t_1, t_2) \in \gamma.$

There are known \cite {Gakh} the following formulas, which are the multidimensional analogues of the Sokhotsky - Plemel formulas:
$$
\Phi^{++} \pm \Phi^{+-} \pm \Phi^{-+} + \Phi^{--} = 
\left \{
\begin{array}{cc}
S_{12} \varphi, \\ {}
\varphi,\\
\end{array}
\right.
\eqno (6.4)
$$
$$
\Phi^{++} \mp \Phi^{+-} \pm \Phi^{-+} - \Phi^{--} = 
\left \{
\begin{array}{cc}
S_1 \varphi, \\ {}
S_2 \varphi.\\
\end{array}
\right.
\eqno (6.5)
$$

Using the definition of the bihypersingular integral, we have
$$
-\frac{1}{\pi^2} \int\limits_{\gamma_1} \int\limits_{\gamma_2} 
\frac{x(\tau_1, \tau_2)d\tau_1 d\tau_2}{(\tau_1-t_1)^p (\tau_2-t_2)^p} =
$$
$$
= -\frac{1}{((p-1)!)^2} \frac{1}{\pi^2} 
\int\limits_{\gamma_1} \int\limits_{\gamma_2} 
\frac{x^{(p-1, p-1)}(\tau_1, \tau_2)d\tau_1 d\tau_2}{(\tau_1-t_1) (\tau_2-t_2)}.
\eqno (6.6)
$$

It follows from (6.4) that
$$
-\frac{1}{\pi^2} \int\limits_{\gamma_1} \int\limits_{\gamma_2} 
\frac{x^{(p-1, p-1)}(\tau_1, \tau_2)d\tau_1 d\tau_2}{(\tau_1-t_1) (\tau_2-t_2)}=
$$
$$
= X^{++}_{p-1,p-1}(t_1,t_2) + X^{+-}_{p-1,p-1}(t_1,t_2) + X^{-+}_{p-1,p-1}(t_1,t_2) + X^{--}_{p-1,p-1}(t_1,t_2),
\eqno (6.7)
$$
where
$$
X_{p-1,p-1}(z_1,z_2) = \frac{1}{(2\pi i)^2} \int\limits_{\gamma_1}\int\limits_{\gamma_2} 
\frac{x^{(p-1, p-1)}(\tau_1, \tau_2)d\tau_1 d\tau_2}{(\tau_1-z_1) (\tau_2-z_2)}.
$$

Since the function $x^{(p-1, p-1)}(z_1,z_2)$ is assumed to be analytic in the domain $ G$, then from the Cauchy integral formula \cite{Fuk} 
follows that
$$
X^{++}_{p-1,p-1}(t_1,t_2) = x^{(p-1, p-1)}(t_1,t_2), 
$$
$$
X^{\pm \mp}_{p-1,p-1}(t_1,t_2) = 0,
$$
$$
X^{--}_{p-1,p-1}(t_1,t_2) = 0.
$$

From this and (6.6), (6.7) we have
$$
-\frac{1}{\pi^2} \int\limits_{\gamma_1} \int\limits_{\gamma_2} 
\frac{x(\tau_1, \tau_2)d\tau_1 d\tau_2}{(\tau_1-t_1)^p (\tau_2-t_2)^p} =
\frac{1}{((p-1)!)^2}x^{(p-1, p-1)}(t_1,t_2).
\eqno (6.8)
$$

It was previously shown that 
$$
\frac{1}{\pi i} \int\limits_{\gamma_1} \frac{x(\tau_1, t_2)}{ (\tau_1-t_1)^p}d\tau_1 =
\frac{1}{(p-1)!}\frac{\partial^{p-1}x(t_1,t_2)}{\partial t_1^{p-1}},
\eqno (6.9)
$$
$$
\frac{1}{\pi i} \int\limits_{\gamma_2} \frac{x(t_1,\tau_2)}{(\tau_2-t_2)^p}d\tau_2 =
\frac{1}{(p-1)!}\frac{\partial^{p-1}x(t_1,t_2)}{\partial t_2^{p-1}}.
\eqno (6.10)
$$

Substituting (6.8) - (6.10) into equation (6.1), we receive the differential equation
$$
a(t_1,t_2) x(t_1,t_2) + \frac{b(t_1,t_2)}{(p-1)!}
\frac{\partial^{p-1}x(t_1,t_2)}{\partial t_1^{p-1}} +
$$
$$
+ \frac{c(t_1,t_2)}{(p-1)!}
\frac{\partial^{p-1}x(t_1,t_2)}{\partial t_2^{p-1}} +
\frac{d(t_1,t_2)}{((p-1)!)^2}
\frac{\partial^{2p-2}x(t_1,t_2)}{\partial t_1^{p-1}\partial t_2^{p-1}} = f(t_1,t_2).
\eqno (6.11)
$$

Thus, it is shown that if the hypersingular integral equation (6.1) has a solution analytic in the domain $ G$, then it is also a solution of the differential equation (6.11). It is easy to see that the converse also holds: if the differential equation (6.11) has a solution analytic in the domain $ G$, then it is also a solution of the hypersingular equation (6.1).

{\bf Example.}

Consider the bihypersingular integral equation
$$
a\int\limits_{\gamma_1} \int\limits_{\gamma_2} 
\frac{x(\tau_1, t_2)}{(\tau_1-t_1)^2}d\tau_1 + 
(bt_1 + c)\int\limits_{\gamma_1} \int\limits_{\gamma_2} 
\frac{x(t_1,\tau_2)}{(\tau_2-t_2)^2}d\tau_2 = 0,
\eqno (6.12)
$$
where $ \gamma_i  $ is a closed smooth contour in the plane of the complex variable $ z_i, $ $ i = 1,2, \ a, b, c  $ are constants.

We denote by $ G = G_1 \times G_2 $ the domain inside which the contour $\gamma = \gamma_1 \times \gamma_2,$ $\gamma_i \in G_i,$ $i=1,2,$ is.

It follows from the above that if equation (6.12) has an analytic solution in the domain $ G$, then this solution is also a solution of the differential equation
$$
a\frac{\partial x(t_1,t_2)}{\partial t_1} + (b t_1 + c) \frac{\partial x(t_1,t_2)}{\partial t_2}=0.
\eqno (6.13)
$$

This equation has the solution \cite{Pol}
$$
X(t_1,t_2)=\frac12 bt_1^2+ct_1-at_2.
$$

One can see that the function $ X (t_1, t_2) $ is also the solution of the hypersingular integral equation (6.12).

Example. We consider the nonlinear bihypersingular integral equation
$$
\frac{1}{\pi i} \int\limits_{\gamma_1}\frac{x(\tau_1, t_2)}{(\tau_1-t_1)^3}d\tau_1
+ \frac{1}{\pi i} \int\limits_{\gamma_2}\frac{x(t_1, \tau_2)}{(\tau_2-t_2)^3}d\tau_2 - kx^n(t_1,t_2) = 0, \, t_i \in \gamma_i, \, i=1,2,
\eqno (6.14)
$$
where $\gamma_i $ is a bounded smooth closed contour on the plane of the complex variable $ z_i, $ $ i = 1,2, $ $ n - $ integer, $ n \geq 2. $

Equation (6.14) is reduced to the partial differential equation of the following form
$$
\frac{\partial^2 x(t_1,t_2)}{\partial t^2_1} + \frac{\partial^2 x(t_1,t_2)}{\partial t^2_2} -
k x^n(t_1,t_2)=0,
\eqno (6.15)
$$
where $t_i \in \gamma_i,$ $i=1,2.$

It is assumed that equation (6.14) has a solution that is an analytic function in the domain $ G $, including the contours $ \gamma_1 $ and $ \gamma_2. $

Under condition that the variables $ t_1 $ and $ t_2 $ are real, equation (6.15) has set of solutions. We will consider exact solutions of the form \cite{Pol}
$$
x(t_1,t_2) = (At_1 + Bt_2 + C)^{2/(1-n)}, \, B = \pm \sqrt{\frac{k(n-1)}{2(n+1)}-A^2};
\eqno (6.16)
$$
$$
x(t_1,t_2) = s[(t_1+C_1)^2 + (t_2+C_2)^2]^{1/(1-n)}, \, s=[\frac{1}{4} k(1-n)^2]^{1/(1-n)},
\eqno (6.17)
$$
where $A, C, C_1, C_2 $ are constants.

Consider the solution (6.16). The choice of the constant $ A $ is bounded by the condition
$0 \leq A \leq \sqrt{\frac{k(n-1)}{2(n+1)}}.$ We choose  constant $ C $ such that $(t_1,t_2) \in G$ $(A t_1 + B t_2 +C) \ne 0.$ (Here it is assumed that $ t_i, $ $ i = 1,2,  $ are complex variables). Then the function $ x (t_1, t_2) $ is analytic in the domain $ G. $ By  direct verification one can see that the function $ x (t_1, t_2) $ satisfies to the equation (6.14).

Similarly, for constants $ C_1 $ and $ C_2 $ such that the function
$ (t_1 + C_1)^2 + (t_2 + C_2)^2 \ne 0 $ for $ (t_1, t_2) \in G,$ the function (6.17) satisfies the equation (6.14).

Summary. The article presents a method for converting hypersingular integral equations into differential equations. The method allows one to obtain solutions of wide classes of hypersingular and polyhypersingular integral equations in closed form.

The work has been supported with Russian Foundation for Basic Research (Grant 16-01-00594).

\end{document}